# Generalized backward doubly stochastic differential equations and SPDEs with nonlinear Neumann boundary conditions

BRAHIM BOUFOUSSI[1,*], JAN VAN CASTEREN[2] and N. MRHARDY[1,**]

[1]*Department of Mathematics, Faculty of Sciences Semlalia, Cadi Ayyad University 2390 Marrakesh, Morocco. E-mail: [*]boufoussi@ucam.ac.ma; [**]n.mrhardy@ucam.ac.ma*
[2]*Department of Mathematics and Computer Science, University of Antwerp, Middelheimlaan 1, 2020 Antwerp, Belgium. E-mail: jan.vancasteren@ua.ac.be*

In this paper a new class of generalized backward doubly stochastic differential equations is investigated. This class involves an integral with respect to an adapted continuous increasing process. A probabilistic representation for viscosity solutions of semi-linear stochastic partial differential equations with a Neumann boundary condition is given.

## 1. Introduction

Backward stochastic differential equations (BSDEs) were introduced by Pardoux and Peng [6], and it was shown in various papers that stochastic differential equations (SDEs) of this type give a probabilistic representation for the solution (at least in the viscosity sense) of a large class of system of semi-linear parabolic partial differential equations (PDEs). A new class of BSDEs, called backward doubly stochastic (BDSDEs), was considered by Pardoux and Peng [7]. This new kind of BSDEs seems to be suitable giving a probabilistic representation for a system of parabolic stochastic partial differential equations (SPDEs). We refer to Pardoux and Peng [7] for the link between SPDEs and BDSDEs in the particular case where solutions of SPDEs are regular. The more general situation is much more delicate to treat because of difficulties of extending the notion of stochastic viscosity solutions to SPDEs.

The notion of viscosity solution for PDEs was introduced by Crandall and Lions [3] for certain first-order Hamilton–Jacobi equations. Today the theory has become an important tool in many applied fields, especially in optimal control theory and numerous subjects related to it.

The stochastic viscosity solution for semi-linear SPDEs was introduced for the first time in Lions and Souganidis [4]. They use the so-called 'stochastic characteristics' to remove the stochastic integrals from an SPDE. Another way of defining a stochastic viscosity solution of SPDEs is via an appeal to the Doss–Sussman transformation. Buckdahn and Ma [2] were the first to use this approach in order to connect the stochastic viscosity solution of SPDEs with BDSDEs. The aim of this paper is to refer to the technique of Buckdahn and Ma [2] to establish the existence result for semi-linear SPDE with Neumann boundary condition of the form:

$$\begin{cases} \mathrm{d}u(t,x) + [Lu(t,x) + f(t,x,u(t,x),\sigma^*(x)\nabla u(t,x))]\,\mathrm{d}t \\ \quad + \sum_{i=1}^{d} g_i(t,x,u(t,x))\,\overleftarrow{\mathrm{d}B_t^i} = 0, \qquad (t,x) \in [0,T] \times \mathrm{G}, \\ u(T,x) = l(x), \qquad\qquad\qquad\qquad\qquad\quad x \in \mathbb{R}^d, \\ \dfrac{\partial u}{\partial n}(t,x) + h(t,x,u(t,x)) = 0, \qquad\qquad x \in \partial\mathrm{G}. \end{cases}$$

Here $B$ is a standard Brownian motion, $L$ is an infinitesimal generator of some diffusion, G is a connected bounded domain, and $f$, $g$, $l$ and $h$ are some measurable functions. More precisely, we give some direct links between the stochastic viscosity solution of the above SPDE and the solution of the following generalized BDSDE:

$$Y_t = \xi + \int_t^T f(s,Y_s,Z_s)\,\mathrm{d}s + \int_t^T h(s,Y_s)\,\mathrm{d}k_s + \int_t^T g(s,Y_s,Z_s)\,\overleftarrow{\mathrm{d}B_s}$$
$$\quad - \int_t^T Z_s\,\mathrm{d}W_s, \qquad 0 \le t \le T,$$

where $\xi$ is the terminal value and $k$ is a real-valued increasing process. Note that our work can be considered as a generalization of the results obtained by Pardoux and Zhang [8], where the authors treat deterministic PDEs with nonlinear Neumann boundary conditions. In light of the approximation result of Boufoussi and Van Casteren [1] for PDEs, a motivation for establishing a connection between SPDEs and BDSDEs is to give a similar (approximation) result for a semi-linear SPDE with a Neumann boundary condition.

The present paper is organized as follows. An existence and uniqueness result for solutions to generalized BDSDEs is shown in Section 2. In Section 3 we introduce the Doss–Sussman transformation which allows us to give a definition of a stochastic viscosity solution to our SPDE. The existence for such a solution via a corresponding BDSDE is given in Section 4.

## 2. Generalized backward doubly stochastic differential equations

### 2.1. Notation and assumptions

Let $T$ be a fixed final time. Throughout this paper $\{W_t,\ 0 \leq t \leq T\}$ and $\{B_t,\ 0 \leq t \leq T\}$ will denote two independent $d$-dimensional Brownian motions ($d \geq 1$), defined on the complete probability spaces $(\Omega_1, \mathcal{F}_1, \mathbb{P}_1)$ and $(\Omega_2, \mathcal{F}_2, \mathbb{P}_2)$, respectively. For any process $(U_s\colon 0 \leq s \leq T)$ defined on $(\Omega_i, \mathcal{F}_i, \mathbb{P}_i)$ ($i = 1, 2$), we write $\mathcal{F}^U_{s,t} := \sigma\{U_r - U_s, s \leq r \leq t\}$ and $\mathcal{F}^U_t := \mathcal{F}^U_{0,t}$. Unless otherwise specified, we consider:

$$\Omega \triangleq \Omega_1 \times \Omega_2, \qquad \mathcal{F} \triangleq \mathcal{F}_1 \otimes \mathcal{F}_2, \qquad \mathbb{P} \triangleq \mathbb{P}_1 \otimes \mathbb{P}_2.$$

In addition, we put

$$\mathcal{F}_t \triangleq \mathcal{F}^W_t \otimes \mathcal{F}^B_{t,T} \vee \mathcal{N},$$

where $\mathcal{N}$ is the collection of $\mathbb{P}$-null sets. In other words, the $\sigma$-fields $\mathcal{F}_t$, $0 \leq t \leq T$, are $\mathbb{P}$-complete. We notice that the family of $\sigma$-algebras $\mathrm{F} = \{\mathcal{F}_t\}_{0 \leq t \leq T}$ is neither increasing nor decreasing; in particular, it is not a filtration. Random variables $\zeta(\omega)$, $\omega \in \Omega_1$, and $\eta(\omega')$, $\omega' \in \Omega_2$, are considered as random variables on $\Omega$ via the following identification:

$$\zeta(\omega, \omega') = \zeta(\omega), \qquad \eta(\omega, \omega') = \eta(\omega').$$

Let $\{k_t,\ 0 \leq t \leq T\}$ be a continuous, increasing and $\mathcal{F}_t$-adapted real-valued process such that $k_0 = 0$. For any $n \geq 1$, we consider the following spaces of processes:

1. The Banach space $\mathcal{M}^2(\mathrm{F}, [0, T]; \mathbb{R}^n)$ of all equivalence classes (with respect to the measure $\mathrm{d}\mathbb{P} \times \mathrm{d}t$) where each equivalence class contains an $n$-dimensional jointly measurable random process $\{\varphi_t,\ t \in [0, T]\}$ which satisfies:
   (i) $\mathbb{E} \int_0^T |\varphi_t|^2 \, \mathrm{d}t < \infty$;
   (ii) $\varphi_t$ is $\mathcal{F}_t$-measurable, for $\mathrm{d}t$-almost all $t \in [0, T]$. Usually an equivalence class will be identified with (one of) its members.
2. The Banach space $\mathcal{K}^2(\mathrm{F}, [0, T]; \mathbb{R}^n)$ of all (equivalence classes of) $n$-dimensional jointly measurable random processes $\{\varphi_t, t \in [0, T]\}$ which satisfy:

    (i) $\mathbb{E} \int_0^T |\varphi_t|^2 \, \mathrm{d}k_t < \infty$;
    (ii) $\varphi_t$ is $\mathcal{F}_t$-measurable, for almost all $t \in [0, T]$.

    Here equivalence is taken with respect to the measure $\mathrm{d}\mathbb{P} \times \mathrm{d}k_t$.
3. The set $\mathcal{S}^2(\mathrm{F}, [0, T]; \mathbb{R}^n)$ of continuous $n$-dimensional random processes which satisfy:
   (i) $\mathbb{E} \sup_{0 \leq t \leq T} |\varphi_t|^2 < \infty$;
   (ii) $\varphi_t$ is $\bar{\mathcal{F}}_t$-measurable, for almost all $t \in [0, T]$.

We consider coefficients $f$, $g$ and $h$ with the following properties:

$$f:\Omega \times [0,T] \times \mathbb{R}^n \times \mathbb{R}^{n\times d} \longrightarrow \mathbb{R}^n,$$
$$g:\Omega \times [0,T] \times \mathbb{R}^n \times \mathbb{R}^{n\times d} \longrightarrow \mathbb{R}^{n\times d},$$
$$h:\Omega \times [0,T] \times \mathbb{R}^n \longrightarrow \mathbb{R}^n,$$

such that there exist $\mathcal{F}_t$-adapted processes $\{f_t, g_t, h_t : 0 \leq t \leq T\}$ with values in $[1, +\infty)$ and with the property that for any $(t, y, z) \in [0,T] \times \mathbb{R}^n \times \mathbb{R}^{n\times d}$, and $\mu > 0$, the following hypotheses are satisfied for some strictly positive finite constant $K$:

$$(\text{H}_1) \quad \begin{cases} f(t,y,z), g(t,y,z) \text{ and } h(t,y) \text{ are } \mathcal{F}_t\text{-measurable processes,} \\ |f(t,y,z)| \leq f_t + K(|y| + \|z\|), \\ |g(t,y,z)| \leq g_t + K(|y| + \|z\|), \\ |h(t,y)| \leq h_t + K|y|, \\ \mathbb{E}\left(\int_0^T e^{\mu k_t} f_t^2 \, dt + \int_0^T e^{\mu k_t} g_t^2 \, dt + \int_0^T e^{\mu k_t} h_t^2 \, dk_t\right) < \infty. \end{cases}$$

Moreover, we assume that there exist constants $c > 0$, $\beta_1 > 0$ and $0 < \alpha < 1$ such that for any $(y_1, z_1), (y_2, z_2) \in \mathbb{R}^n \times \mathbb{R}^{n\times d}$,

$$(\text{H}_2) \quad \begin{cases} (\text{i}) & |f(t,y_1,z_1) - f(t,y_2,z_2)|^2 \leq c(|y_1-y_2|^2 + \|z_1-z_2\|^2), \\ (\text{ii}) & |g(t,y_1,z_1) - g(t,y_2,z_2)|^2 \leq c|y_1-y_2|^2 + \alpha \|z_1-z_2\|^2, \\ (\text{iii}) & |h(t,y_1) - h(t,y_2)| \leq \beta_1 |y_1-y_2|. \end{cases}$$

Throughout this paper, $\langle \cdot, \cdot \rangle$ will denote the scalar product on $\mathbb{R}^n$, i.e. $\langle x, y \rangle := \sum_{i=1}^n x_i y_i$, for all $(x,y) \in \mathbb{R}^n \times \mathbb{R}^n$. Sometimes, we will also use the notation $x^\star y$ to designate $\langle x, y \rangle$.

Finally, $C$ will always denote a finite constant whose value may change from one line to the next, and which usually is (strictly) positive.

## 2.2. Existence and uniqueness theorem

Suppose that we are given a terminal condition $\xi \in L^2(\Omega, \mathcal{F}_T, \mathbb{P})$ such that, for all $\mu > 0$,

$$\mathbb{E}(e^{\mu k_T} |\xi|^2) < \infty.$$

**Definition 2.1.** *By definition, a solution to a generalized BDSDE $(\xi, f, g, h, k)$ is a pair $(Y, Z) \in \mathcal{S}^2(\mathrm{F}, [0,T]; \mathbb{R}^n) \times \mathcal{M}^2(\mathrm{F}, [0,T]; \mathbb{R}^{n\times d})$, such that, for any $0 \leq t \leq T$,*

$$Y_t = \xi + \int_t^T f(s, Y_s, Z_s) \, ds + \int_t^T h(s, Y_s) \, dk_s + \int_t^T g(s, Y_s, Z_s) \overleftarrow{dB_s} - \int_t^T Z_s \, dW_s. \quad (1)$$

*Here $\overleftarrow{dB_s}$ denotes the classical backward Itô integral with respect to the Brownian motion $B$.*

**Remark 2.1.** If $h$ satisfies $(H_2)(iii)$ then, by changing the solutions and the coefficients $f$, $g$ and $h$, we may and do suppose that $h$ satisfies a stronger condition of the form

(iv) $\langle y_1 - y_2, h(t, y_1) - h(t, y_2) \rangle \leq \beta_2 |y_1 - y_2|^2$, where $\beta_2 < 0$.

Indeed, $(Y_t, Z_t)$ solves the generalized BDSDE in (1) if and only if for every (some) $\eta > 0$ the pair $(\overline{Y}_t, \overline{Z}_t) = (e^{\eta k_t} Y_t, e^{\eta k_t} Z_t)$ solves an analogous generalized BDSDE, with $f$, $g$ and $h$ replaced respectively by

$$\overline{f}(t, y, z) = e^{\eta k_t} f(t, e^{-\eta k_t} y, e^{-\eta k_t} z);$$
$$\overline{g}(t, y, z) = e^{\eta k_t} g(t, e^{-\eta k_t} y, e^{-\eta k_t} z);$$
$$\overline{h}(t, y) = e^{\eta k_t} h(t, e^{-\eta k_t} y) - \eta y.$$

Then we can always choose $\eta$ such that the function $\overline{h}$ satisfies (iv) with a strictly negative $\beta_2$.

Our main goal in this section is to prove the following theorem.

**Theorem 2.1.** *Under the above hypotheses* $(H_1)$ *and* $(H_2)$ *there exists a unique solution for the generalized BDSDE in* (1).

We will follow the same line of arguments as Pardoux and Peng [7] did. So let us first establish the result in Theorem 2.1 for BDSDEs where the coefficients $f$, $g$ and $h$ do not depend on $(y, z)$. More precisely, let $f$, $h : \Omega \times [0, T] \to \mathbb{R}^n$ and $g : \Omega \times [0, T] \to \mathbb{R}^{n \times d}$ satisfy $(H_1)$, and let $\xi$ and $k$ be as before. Consider the equation:

$$Y_t = \xi + \int_t^T f(s)\,\mathrm{d}s + \int_t^T h(s)\,\mathrm{d}k_s + \int_t^T g(s)\,\overleftarrow{\mathrm{d}B_s} - \int_t^T Z_s\,\mathrm{d}W_s. \tag{2}$$

Then we have the following result.

**Theorem 2.2.** *Under hypothesis* $(H_1)$, *there exists a unique solution to equation* (2).

**Proof.** To show the existence, we consider the filtration $\mathcal{G}_t = \mathcal{F}_t^W \otimes \mathcal{F}_T^B$ and the martingale

$$M_t = \mathbb{E}\left[\xi + \int_0^T f(s)\,\mathrm{d}s + \int_0^T h(s)\,\mathrm{d}k_s + \int_0^T g(s)\,\overleftarrow{\mathrm{d}B_s} / \mathcal{G}_t\right], \tag{3}$$

which is clearly a square integrable martingale by $(H_1)$. As in Pardoux and Peng [7], an extension of Itô's martingale representation theorem yields the existence of a $\mathcal{G}_t$-progressively measurable process $(Z_t)$ with values in $\mathbb{R}^{n \times d}$ such that

$$\mathbb{E}\int_0^T \|Z_t\|^2\,\mathrm{d}t < \infty \quad \text{and} \quad M_T = M_t + \int_t^T Z_s\,\mathrm{d}W_s. \tag{4}$$

We subtract the quantity $\int_0^t f(s)\,\mathrm{d}s + \int_0^t h(s)\,\mathrm{d}k_s + \int_0^t g(s)\,\overleftarrow{\mathrm{d}B_s}$ from both sides of the martingale in (3) and employ the martingale representation in (4) to obtain

$$Y_t = \xi + \int_t^T f(s)\,\mathrm{d}s + \int_t^T h(s)\,\mathrm{d}k_s + \int_t^T g(s)\,\overleftarrow{\mathrm{d}B_s} - \int_t^T Z_s\,\mathrm{d}W_s,$$

where

$$Y_t = \mathbb{E}\left[\xi + \int_t^T f(s)\,\mathrm{d}s + \int_t^T h(s)\,\mathrm{d}k_s + \int_t^T g(s)\,\overleftarrow{\mathrm{d}B_s}/\mathcal{G}_t\right].$$

It remains to prove the uniqueness and to show that $Y_t$ and $Z_t$ are $\mathcal{F}_t$-measurable; the proof is analogous to that of Pardoux and Peng [7]: see Proposition 1.3, and is therefore omitted. □

We will also need the following generalized Itô formula. In the proof we use arguments which are similar to those used by Pardoux and Peng in [7].

**Lemma 2.3.** *Let* $\alpha \in \mathcal{S}^2(\mathrm{F},[0,T];\mathbb{R}^n)$, $\beta \in \mathcal{M}^2(\mathrm{F},[0,T];\mathbb{R}^n)$, $\gamma \in \mathcal{M}^2(\mathrm{F},[0,T];\mathbb{R}^{n\times d})$, $\theta \in \mathcal{K}^2(\mathrm{F},[0,T];\mathbb{R}^n)$ *and* $\delta \in \mathcal{M}^2(\mathrm{F},[0,T];\mathbb{R}^{n\times d})$ *be such that*

$$\alpha_t = \alpha_0 + \int_0^t \beta_s\,\mathrm{d}s + \int_0^t \theta_s\,\mathrm{d}k_s + \int_0^t \gamma_s\,\overleftarrow{\mathrm{d}B_s} + \int_0^t \delta_s\,\mathrm{d}W_s.$$

*Then, for any function* $\phi \in \mathcal{C}^2(\mathbb{R}^n,\mathbb{R})$,

$$\phi(\alpha_t) = \phi(\alpha_0) + \int_0^t \langle\nabla\phi(\alpha_s),\beta_s\rangle\,\mathrm{d}s + \int_0^t \langle\nabla\phi(\alpha_s),\theta_s\rangle\,\mathrm{d}k_s + \int_0^t \langle\nabla\phi(\alpha_s),\gamma_s\,\overleftarrow{\mathrm{d}B_s}\rangle$$
$$+ \int_0^t \langle\nabla\phi(\alpha_s),\delta_s\,\mathrm{d}W_s\rangle - \frac{1}{2}\int_0^t \mathrm{tr}[\phi''(\alpha_s)\gamma_s\gamma_s^*]\,\mathrm{d}s + \frac{1}{2}\int_0^t \mathrm{Tr}[\phi''(\alpha_s)\delta_s\delta_s^*]\,\mathrm{d}s.$$

*In particular,*

$$|\alpha|_t^2 = |\alpha_0|^2 + 2\int_0^t \langle\alpha_s,\beta_s\rangle\,\mathrm{d}s + 2\int_0^t \langle\alpha_s,\theta_s\rangle\,\mathrm{d}k_s + 2\int_0^t \langle\alpha_s,\gamma_s\,\overleftarrow{\mathrm{d}B_s}\rangle$$
$$+ 2\int_0^t \langle\alpha_s,\delta_s\,\mathrm{d}W_s\rangle - \int_0^t \|\gamma_s\|^2\,\mathrm{d}s + \int_0^t \|\delta_s\|^2\,\mathrm{d}s.$$

Next, we establish an a priori estimate for the solution of the BSDE in (1).

**Proposition 2.4.** *Let the conditions* $(\mathrm{H}_1)$ *and* $(\mathrm{H}_2)$ *be satisfied. If* $\{(Y_t,Z_t); 0\leq t\leq T\}$ *is a solution of BDSDE* (1), *then there exists a finite constant* $C$, *which depends on* $K$, $T$ *and* $\beta_2$, *such that for all* $\mu\in\mathbb{R}$ *and* $\lambda>0$ *the following inequality holds:*

$$\mathbb{E}\left(\sup_{0\leq t\leq T} \mathrm{e}^{\mu t+\lambda k_t}|Y_t|^2 + \int_0^T \mathrm{e}^{\mu t+\lambda k_t}|Y_t|^2\,\mathrm{d}k_t + \int_0^T \mathrm{e}^{\mu t+\lambda k_t}\|Z_t\|^2\,\mathrm{d}t\right)$$

$$\leq C\mathbb{E}\left(e^{\mu T+\lambda k_T}|\xi|^2 + \int_0^T e^{\mu t+\lambda k_t}|f_t|^2\,dt + \int_0^T e^{\mu t+\lambda k_t}|h_t|^2\,dk_t + \int_0^T e^{\mu t+\lambda k_t}|g_t|^2\,dt\right).$$

**Proof.** Classical arguments, such as Doob's inequality, justify the fact that the processes $\int_0^t e^{\mu s+\lambda k_s}\langle Y_s, g(s,Y_s,Z_s)\overleftarrow{dB_s}\rangle$ and $\int_0^t e^{\mu s+\lambda k_s}\langle Y_s, Z_s\,dW_s\rangle$ are uniformly integrable martingales. By Lemma 2.3 we then have

$$\mathbb{E}\left[e^{\mu t+\lambda k_t}|Y_t|^2 + \int_t^T e^{\mu s+\lambda k_s}\|Z_s\|^2\,ds + \lambda\int_t^T e^{\mu s+\lambda k_s}|Y_s|^2\,dk_s\right]$$

$$\leq \mathbb{E}\bigg[e^{\mu T+\lambda k_T}|\xi|^2 + 2\int_t^T e^{\mu s+\lambda k_s}\langle Y_s, f(s,Y_s,Z_s)\rangle\,ds$$

$$+ 2\int_t^T e^{\mu s+\lambda k_s}\langle Y_s, h(s,Y_s)\rangle\,dk_s + \int_t^T e^{\mu s+\lambda k_s}|g(s,Y_s,Z_s)|^2\,ds$$

$$- \mu\int_t^T e^{\mu s+\lambda k_s}|Y_s|^2\,ds\bigg]. \tag{5}$$

But from (H$_1$), (H$_2$) and the fact that

$$2ab \leq \frac{1-\alpha}{2c}a^2 + \frac{2c}{1-\alpha}b^2, \qquad c>0,$$

it follows that there exists a constant $c(\alpha)$ such that

$$2\langle y, f(s,y,z)\rangle \leq c|f_s|^2 + c(\alpha)|y|^2 + \frac{1-\alpha}{2}\|z\|^2, \tag{6}$$

$$2\langle y, h(s,y)\rangle \leq 2\beta_2|y|^2 + |y|\times|h_s| \leq (2\beta_2+|\beta_2|)|y|^2 + \frac{1}{|\beta_2|}h_s^2, \tag{7}$$

$$\|g(s,y,z)\|^2 \leq c|y|^2 + \alpha\|z\|^2 + (1+\varepsilon)g_s^2 + \frac{c}{\varepsilon}|y|^2 + \frac{\alpha}{\varepsilon}\|z\|^2. \tag{8}$$

Inserting $\varepsilon = 3\alpha/(1-\alpha)$ into (8) replaces the latter inequality by

$$\|g(s,y,z)\|^2 \leq c'(\alpha)|y|^2 + \left(\alpha+\frac{1-\alpha}{3}\right)\|z\|^2 + \frac{1+2\alpha}{1-\alpha}g_s^2.$$

Consequently, by (5) we obtain for the same constant $c(\alpha)$ the inequality

$$\mathbb{E}\left(e^{\mu t+\lambda k_t}|Y_t|^2 + (\lambda+|\beta_2|)\int_t^T e^{\mu s+\lambda k_s}|Y_s|^2\,dk_s + \frac{1-\alpha}{6}\int_t^T e^{\mu s+\lambda k_s}\|Z_s\|^2\,ds\right)$$

$$\leq \mathbb{E}\bigg(e^{\mu T+\lambda k_T}|\xi|^2 + (c(\alpha)-\mu)\int_t^T e^{\mu s+\lambda k_s}|Y_s|^2\,ds + c\int_t^T e^{\mu s+\lambda k_s}|f_s|^2\,ds$$

$$+ \frac{1}{c}\int_t^T e^{\mu s+\lambda k_s}|h_s|^2\,dk_s + \frac{1+2\alpha}{1-\alpha}\int_t^T e^{\mu s+\lambda k_s}|g_s|^2\,ds\bigg).$$

Then, from Gronwall's lemma we obtain

$$\sup_{0 \le t \le T} \mathbb{E}\left(e^{\mu t + \lambda k_t}|Y_t|^2 + \int_0^T e^{\mu s + \lambda k_s}|Y_s|^2 \, dk_s + \int_0^T e^{\mu s + \lambda k_s}\|Z_s\|^2 \, ds\right)$$

$$\le C\mathbb{E}\left(e^{\mu T + \lambda k_T}|\xi|^2 + \int_0^T e^{\mu s + \lambda k_s}|f_s|^2 \, ds \right.$$

$$\left. + \int_0^T e^{\mu s + \lambda k_s}|h_s|^2 \, dk_s + \int_0^T e^{\mu s + \lambda k_s}|g_s|^2 \, ds\right). \tag{9}$$

Finally, Proposition 2.4 follows from the Burkholder–Davis–Gundy inequality and (9). □

Next, let $(\xi, f, g, h, k)$ and $(\xi', f', g', h', k')$ be two sets of data, each satisfying conditions $(H_1)$ and $(H_2)$. Then we have the following result:

**Proposition 2.5.** *Let $(Y, Z)$ (or $(Y', Z')$) denote a solution of the BDSDE$(\xi, f, g, h, k)$ (or BDSDE$(\xi', f', g', h', k')$). With the notation*

$$(\overline{Y}, \overline{Z}, \overline{\xi}, \overline{f}, \overline{g}, \overline{h}, \overline{k}) = (Y - Y', Z - Z', \xi - \xi', f - f', g - g', h - h', k - k'),$$

*it follows that for every $\mu > 0$, there exists a constant $C > 0$ such that*

$$\mathbb{E}\left(\sup_{0 \le t \le T} e^{\mu A_t}|\overline{Y}_t|^2 + \int_0^T e^{\mu A_t}\|\overline{Z}_t\|^2 \, dt\right)$$

$$\le C\mathbb{E}\left(e^{\mu A_T}|\overline{\xi}|^2 + \int_0^T e^{\mu A_t}|f(t, Y_t, Z_t) - f'(t, Y_t, Z_t)|^2 \, dt + \int_0^T e^{\mu A_t}|h(t, Y_t)|^2 \, d|\overline{k}|_t \right.$$

$$\left. + \int_0^T e^{\mu A_t}|h(t, Y_t) - h'(t, Y_t)|^2 \, dk'_t + \int_0^T e^{\mu A_t}\|g(t, Y_t, Z_t) - g'(t, Y_t, Z_t)\|^2 \, dt\right).$$

*Here $A_t \triangleq |\overline{k}|_t + k'_t$ and $|\overline{k}|_t$ is the total variation of the process $\overline{k}$.*

**Proof.** The proof follows the same ideas and arguments as in Pardoux and Zhang [8], Proposition 1.2, so we just repeat the main steps. From Lemma 2.3 we obtain

$$e^{\mu A_t}|\overline{Y}_t|^2 + \int_t^T e^{\mu A_s}\|\overline{Z}_s\|^2 \, ds + \mu \int_t^T e^{\mu A_s}|\overline{Y}_s|^2 \, dA_s$$

$$= e^{\mu A_T}|\overline{\xi}|^2 + 2\int_t^T e^{\mu A_s}\langle \overline{Y}_s, f(s, Y_s, Z_s) - f'(s, Y'_s, Z'_s)\rangle \, ds + 2\int_t^T e^{\mu A_s}\langle \overline{Y}_s, h(s, Y_s)\rangle \, d\overline{k}_s$$

$$+ 2\int_t^T e^{\mu A_s}\langle \overline{Y}_s, h(s, Y_s) - h'(s, Y'_s)\rangle \, dk'_s + \int_t^T e^{\mu A_s}\|g(s, Y_s, Z_s) - g'(s, Y'_s, Z'_s)\|^2 \, ds$$

$$+ 2\int_t^T e^{\mu A_s}\langle \overline{Y}_s, g(s, Y_s, Z_s) - g'(s, Y'_s, Z'_s)\rangle \, \overleftarrow{dB}_s - 2\int_t^T e^{\mu A_s}\langle \overline{Y}_s, \overline{Z}_s \, dW_s\rangle. \tag{10}$$

Using conditions (H$_1$), (H$_2$), and the algebraic inequality $2ab \leq a^2/\varepsilon + \varepsilon b^2$, then from (10) we obtain

$$\mathbb{E}\left(e^{\mu A_t}|\overline{Y}_t|^2 + \int_t^T e^{\mu A_s}\|\overline{Z}_s\|^2 \, ds + \mu \int_t^T e^{\mu A_s}|\overline{Y}_s|^2 \, dA_s\right)$$

$$\leq \mathbb{E}\left(e^{\mu A_T}|\overline{\xi}|^2 + C(\alpha)\int_t^T e^{\mu A_s}|\overline{Y}_s|^2 \, ds + \int_t^T e^{\mu A_s}|\overline{f}(s, Y_s, Z_s)|^2 \, ds\right.$$

$$+ \frac{1}{\varepsilon}\int_t^T e^{\mu A_s}|\overline{h}(s, Y_s)|^2 \, dk'_s + \int_t^T e^{\mu A_s}\|\overline{g}(s, Y_s, Z_s)\|^2 \, ds$$

$$+ \frac{1}{\mu}\int_t^T e^{\mu A_s}|h(s, Y_s)|^2 \, d|\overline{k}|_s$$

$$\left.+ \mu \int_t^T e^{\mu A_s}|\overline{Y}_s|^2 \, d|\overline{k}|_s + (2\beta_2 + \varepsilon)\int_t^T e^{\mu A_s}|\overline{Y}_s|^2 \, dk'_s\right). \quad (11)$$

By choosing $\varepsilon = \mu + 2|\beta_2|$, and using Gronwall's lemma, from (11) we infer that

$$\mathbb{E}\left(e^{\mu A_t}|\overline{Y}_t|^2 + \int_0^T e^{\mu A_t}\|\overline{Z}_t\|^2 \, dt\right)$$

$$\leq C(\alpha, \mu)\mathbb{E}\left(e^{\mu A_T}|\overline{\xi}|^2 + \int_0^T e^{\mu A_s}|\overline{f}(s, Y_s, Z_s)|^2 \, ds\right.$$

$$+ \int_0^T e^{\mu A_s}\|\overline{g}(s, Y_s, Z_s)\|^2 \, ds$$

$$\left.+ \int_0^T e^{\mu A_s}|h(s, Y_s)|^2 \, d|\overline{k}|_s + \int_0^T e^{\mu A_s}|\overline{h}(s, Y_s)|^2 \, dk'_s\right). \quad (12)$$

The proposition follows from (12) and the Burkholder–Davis–Gundy inequality. □

**Remark 2.2.** If we denote by $\mathbb{E}^{\mathcal{F}_t}$ the conditional expectation with respect to $\mathcal{F}_t$, then we can show that for every $\mu, \lambda > 0$, there exists a constant $C > 0$ such that $\forall t \in [0, T]$

$$e^{\mu A_t + \lambda t}|\overline{Y}_t|^2 = \mathbb{E}^{\mathcal{F}_t}(e^{\mu A_t + \lambda t}|\overline{Y}_t|^2)$$

$$\leq C\mathbb{E}^{\mathcal{F}_t}\left(e^{\mu A_T + \lambda T}|\overline{\xi}|^2 + \int_0^T e^{\mu A_s + \lambda s}|\overline{f}(s, Y_s, Z_s)|^2 \, ds\right.$$

$$+ \int_0^T e^{\mu A_s + \lambda s}|\overline{h}(s, Y_s)|^2 \, dk'_s$$

$$+ \int_0^T e^{\mu A_s + \lambda s}|h(s, Y_s)|^2 \, d|\overline{k}|_s$$

$$+ \int_0^T e^{\mu A_s + \lambda s} \|\overline{g}(s, Y_s, Z_s)\|^2 \, ds \bigg), \qquad \mathbb{P}\text{-almost surely.}$$

**Proof of Theorem 2.1.** The uniqueness is a consequence of Proposition 2.5. We now turn to the existence. In the space $\mathcal{S}^2(F, [0, T]; \mathbb{R}^n) \times \mathcal{M}^2(F, [0, T]; \mathbb{R}^n)$ we define by recursion the sequence $\{(Y_t^i, Z_t^i)\}_{i=0,1,2,\ldots}$ as follows. Put $Y_t^0 = 0$, $Z_t^0 = 0$. Given the pair $(Y_t^i, Z_t^i)$, we define $f^{i+1}(s) = f(s, Y_s^i, Z_s^i)$, $h^{i+1}(s) = h(s, Y_s^i)$ and $g^{i+1}(s) = g(s, Y_s^i, Z_s^i)$. Now, applying (H$_1$), we obtain

$$|h^{i+1}(s)| \leq h_s + K|Y_s^i| \triangleq h_s^{i+1},$$

and by using Proposition 2.4, we obtain

$$\mathbb{E} \int_0^T e^{\mu k_s} (h_s^{i+1})^2 \, dk_s \leq C\mathbb{E} \bigg( \int_0^T e^{\mu k_s} h_s^2 \, dk_s + \int_0^T e^{\mu k_s} |Y_s^i|^2 \, dk_s \bigg) < \infty.$$

By the same arguments one can show that $f^{i+1}$ and $g^{i+1}$ also satisfy (H$_1$). Using Theorem 2.2, we consider the process $\{(Y_t^{i+1}, Z_t^{i+1})\}$ as being the unique solution to the equation

$$Y_t^{i+1} = \xi + \int_t^T f(s, Y_s^i, Z_s^i) \, ds + \int_t^T h(s, Y_s^i) \, dk_s$$
$$+ \int_t^T g(s, Y_s^i, Z_s^i) \, \overleftarrow{dB_s} - \int_t^T Z_s^{i+1} \, dW_s. \qquad (13)$$

We will show that the sequence $\{(Y_t^i, Z_t^i)\}$ converges in the space $\mathcal{S}^2(F, [0, T]; \mathbb{R}^n) \times \mathcal{M}^2(F, [0, T]; \mathbb{R}^n)$ to a pair of processes $(Y_t, Z_t)$ which will be our solution. Indeed, let

$$\overline{Y}_t^{i+1} \triangleq Y_t^{i+1} - Y_t^i, \qquad \overline{Z}_t^{i+1} \triangleq Z_t^{i+1} - Z_t^i.$$

Let $\mu > 0$, $\lambda > 0$. Using Itô's formula, we obtain

$$e^{\mu t + \lambda k_t} |\overline{Y}_t^{i+1}|^2 + \int_t^T e^{\mu s + \lambda k_s} \|\overline{Z}_s^{i+1}\|^2 \, ds$$
$$= 2 \int_t^T e^{\mu s + \lambda k_s} \langle \overline{Y}_s^{i+1}, f(s, Y_s^i, Z_s^i) - f(s, Y_s^{i-1}, Z_s^{i-1}) \rangle \, ds - \mu \int_t^T e^{\mu s + \lambda k_s} |\overline{Y}_s^{i+1}|^2 \, ds$$
$$+ 2 \int_t^T e^{\mu s + \lambda k_s} \langle \overline{Y}_s^{i+1}, h(s, Y_s^i) - h(s, Y_s^{i-1}) \rangle \, dk_s - \mu \int_t^T e^{\mu s + \lambda k_s} |\overline{Y}_s^{i+1}|^2 \, dk_s$$
$$+ \int_t^T e^{\mu s + \lambda k_s} \|g(s, Y_s^i, Z_s^i) - g(s, Y_s^{i-1}, Z_s^{i-1})\|^2 \, ds - 2 \int_t^T e^{\mu s + \lambda k_s} \langle \overline{Y}_s^{i+1}, \overline{Z}_s^{i+1} \rangle \overleftarrow{dB_s}$$
$$+ 2 \int_t^T e^{\mu s + \lambda k_s} \langle \overline{Y}_s^{i+1}, (g(s, Y_s^i, Z_s^i) - g(s, Y_s^{i-1}, Z_s^{i-1})) \, dW_s \rangle.$$

Taking the expectation, we get

$$\mathbb{E} e^{\mu t+\lambda k_t}|\overline{Y}_t^{i+1}|^2 + \mathbb{E}\int_t^T e^{\mu s+\lambda k_s}\|\overline{Z}_s^{i+1}\|^2\,\mathrm{d}s$$

$$= 2\mathbb{E}\int_t^T e^{\mu s+\lambda k_s}\langle \overline{Y}_s^{i+1}, f(s,Y_s^i,Z_s^i) - f(s,Y_s^{i-1},Z_s^{i-1})\rangle\,\mathrm{d}s$$

$$+ 2\mathbb{E}\int_t^T e^{\mu s+\lambda k_s}\langle \overline{Y}_s^{i+1}, h(s,Y_s^i) - h(s,Y_s^{i-1})\rangle\,\mathrm{d}k_s - \mu\mathbb{E}\int_t^T e^{\mu s+\lambda k_s}|\overline{Y}_s^{i+1}|^2\,\mathrm{d}k_s$$

$$+ \mathbb{E}\int_t^T e^{\mu s+\lambda k_s}\|g(s,Y_s^i,Z_s^i) - g(s,Y_s^{i-1},Z_s^{i-1})\|^2\,\mathrm{d}s - \mu\mathbb{E}\int_t^T e^{\mu s+\lambda k_s}|\overline{Y}_s^{i+1}|^2\,\mathrm{d}s.$$

With the same arguments as in the proof of Proposition 2.5, one can show that there exist constants $c_1(\alpha)$, $c_2(\alpha)$, and $\overline{c} > 0$ such that

$$\mathbb{E} e^{\mu t+\lambda k_t}|\overline{Y}_t^{i+1}|^2 + \mathbb{E}\int_t^T e^{\mu s+\lambda k_s}\|\overline{Z}_s^{i+1}\|^2\,\mathrm{d}s$$

$$+ \mathbb{E}\int_t^T e^{\mu s+\lambda k_s}((\mu - c_1(\alpha))|\overline{Y}_s^{i+1}|^2\,\mathrm{d}s + (\lambda - c_2(\alpha))|\overline{Y}_s^{i+1}|^2\,\mathrm{d}k_s)$$

$$\leq \frac{1+\alpha}{2}\left(\overline{c}\mathbb{E}\int_t^T e^{\mu s+\lambda k_s}|\overline{Y}_s^i|^2\,\mathrm{d}s + c\mathbb{E}\int_t^T e^{\mu s+\lambda k_s}|\overline{Y}_s^i|^2\,\mathrm{d}k_s + \mathbb{E}\int_t^T e^{\mu s+\lambda k_s}\|\overline{Z}_s^i\|^2\,\mathrm{d}s\right).$$

Next, we choose $\mu$ and $\lambda$ in such a way that $\mu - c_1(\alpha) = \overline{c}$ and $\lambda - c_2(\alpha) = c$ to obtain

$$\mathbb{E} e^{\mu t+\lambda k_t}|\overline{Y}_t^{i+1}|^2 + \mathbb{E}\int_t^T e^{\mu s+\lambda k_s}\|\overline{Z}_s^{i+1}\|^2\,\mathrm{d}s$$

$$+ c\mathbb{E}\int_t^T e^{\mu s+\lambda k_s}|\overline{Y}_s^{i+1}|^2\,\mathrm{d}k_s + \overline{c}\mathbb{E}\int_t^T e^{\mu s+\lambda k_s}|\overline{Y}_s^{i+1}|^2\,\mathrm{d}s$$

$$\leq \left(\frac{1+\alpha}{2}\right)^i\left[\overline{c}\mathbb{E}\int_t^T e^{\mu s+\lambda k_s}|\overline{Y}_s^1|^2\,\mathrm{d}s + c\mathbb{E}\int_t^T e^{\mu s+\lambda k_s}|\overline{Y}_s^1|^2\,\mathrm{d}k_s + \mathbb{E}\int_t^T e^{\mu s+\lambda k_s}\|\overline{Z}_s^1\|^2\,\mathrm{d}s\right].$$

Since $(1+\alpha)/2 < 1$, then $\{(Y_t^i, Z_t^i)\}_{i=1,\ldots}$ is a Cauchy sequence in the space

$$L^2(\mathrm{F},[0,T];\mathbb{R}^n) \times \mathcal{M}^2(\mathrm{F},[0,T];\mathbb{R}^{n\times d}).$$

From the Burkholder–Davis–Gundy inequality it follows that the sequence $(Y_t^i)$ is also a Cauchy sequence in the space $\mathcal{S}^2(\mathrm{F},[0,T];\mathbb{R}^n)$. By completeness its limit $(Y_t, Z_t) = \lim_{i\to\infty}(Y_t^i, Z_t^i)$ exists in the space $\mathcal{S}^2(\mathrm{F},[0,T];\mathbb{R}^n) \times \mathcal{M}^2(\mathrm{F},[0,T];\mathbb{R}^{n\times d})$. Passing to the limit in equation (13), we obtain the result in Theorem 2.2. □

## 3. Viscosity solutions

In this section we introduce the notion of stochastic viscosity solutions to semi-linear SPDEs with Neumann boundary conditions, and by using the generalized BDSDE we prove the existence of such solutions.

### 3.1. Preliminaries and definitions

With the same notation as in Section 2, let $\mathbf{F}^B \triangleq \{\mathcal{F}^B_{t,T}\}_{0\leq t\leq T}$. By $\mathcal{M}^B_{0,T}$ we will denote all the $\mathbf{F}^B$-stopping times $\tau$ such that $0 \leq \tau \leq T$, $\mathbb{P}$-almost surely. For generic Euclidean spaces $E$ and $E_1$ we introduce the following vector spaces of functions:

1. $\mathcal{C}^{k,\ell}([0,T] \times E; E_1)$ stands for the space of all $E_1$-valued functions defined on $[0,T] \times E$ which are $k$ times continuously differentiable in $t$ and $\ell$ times continuously differentiable in $x$, and $\mathcal{C}^{k,\ell}_b([0,T] \times E; E_1)$ denotes the subspace of $\mathcal{C}^{k,\ell}([0,T] \times E; E_1)$ in which all functions have uniformly bounded partial derivatives.
2. For any sub-$\sigma$-field $\mathcal{G} \subseteq \mathcal{F}^B_T$, $\mathcal{C}^{k,\ell}(\mathcal{G}, [0,T] \times E; E_1)$ (or $\mathcal{C}^{k,\ell}_b(\mathcal{G}, [0,T] \times E; E_1)$) denotes the space of all $\mathcal{C}^{k,\ell}([0,T] \times E; E_1)$ (or $\mathcal{C}^{k,\ell}_b([0,T] \times E; E_1)$-valued) random variables that are $\mathcal{G} \otimes \mathcal{B}([0,T] \times E)$-measurable.
3. $\mathcal{C}^{k,\ell}(\mathbf{F}^B, [0,T] \times E; E_1)$ (or $\mathcal{C}^{k,\ell}_b(\mathbf{F}^B, [0,T] \times E; E_1)$) is the space of all random fields $\alpha \in \mathcal{C}^{k,\ell}(\mathcal{F}^B_T, [0,T] \times E; E_1)$ (or $\mathcal{C}^{k,\ell}_b(\mathcal{F}^B_T, [0,T] \times E; E_1)$), such that for fixed $x \in E$, the mapping $(t,w) \to \alpha(t,\omega,x)$ is $\mathbf{F}^B$-progressively measurable.
4. For any sub-$\sigma$-field $\mathcal{G} \subseteq \mathcal{F}^B_T$ and real number $p \geq 0$, $L^p(\mathcal{G}; E)$ stands for all $E$-valued $\mathcal{G}$-measurable random variables $\xi$ such that $\mathbb{E}|\xi|^p < \infty$.

Furthermore, for $(t,x,y) \in [0,T] \times \mathbb{R}^n \times \mathbb{R}$, we write

$$D_x = \left(\frac{\partial}{\partial x_1}, \ldots, \frac{\partial}{\partial x_n}\right), \qquad D_{xx} = (\partial^2_{x_i x_j})^n_{i,j=1}, \qquad D_y = \frac{\partial}{\partial y}, \qquad D_t = \frac{\partial}{\partial t}.$$

The meaning of $D_{xy}$ and $D_{yy}$ is then self-explanatory.

Let G be an open connected bounded domain of $\mathbb{R}^n$ ($n \geq 1$). We suppose that G is a smooth domain, which is such that for a function $\phi \in \mathcal{C}^2_b(\mathbb{R}^n)$, G and its boundary $\partial G$ are characterized by $G = \{\phi > 0\}$, $\partial G = \{\phi = 0\}$ and, for any $x \in \partial G$, $\nabla \phi(x)$ is the unit normal vector pointing towards the interior of G.

In this section, we consider continuous coefficients $f$ and $h$,

$$f : \Omega_2 \times [0,T] \times \overline{G} \times \mathbb{R} \times \mathbb{R}^d \longrightarrow \mathbb{R},$$
$$h : \Omega_2 \times [0,T] \times \overline{G} \times \mathbb{R} \longrightarrow \mathbb{R}$$

with the property that for all $x \in \overline{G}$, $f(\cdot,x,\cdot,\cdot)$ and $h(\cdot,x,\cdot)$ are Lipschitz continuous in $x$ and satisfy the conditions (H$'_1$) and (H$_2$), uniformly in $x$, where, for some constant

$K > 0$, the condition (H$'_1$) is:

(H$'_1$) $\begin{cases} |f(t,x,y,z)| \leq K(1+|y|+|x|+\|z\|), \\ |h(t,x,y)| \leq K(1+|y|+|x|). \end{cases}$

Furthermore, we shall make use of the following assumptions:

(H$_3$) The functions $\sigma : \mathbb{R}^n \to \mathbb{R}^{n \times d}$ and $b : \mathbb{R}^n \to \mathbb{R}^n$ are uniformly Lipschitz continuous, with common Lipschitz constant $K > 0$.

(H$_4$) The function $l : \overline{G} \to \mathbb{R}$ is continuous such that, for some constant $K > 0$,

$$|l(x)| \leq K(1+|x|), \qquad x \in \overline{G}.$$

(H$_5$) The function $g \in \mathcal{C}_b^{0,2,3}([0,T] \times \overline{G} \times \mathbb{R}; \mathbb{R}^d)$.

We consider the second-order differential operator

$$L = \frac{1}{2} \sum_{i,j=1}^n (\sigma(x)\sigma^*(x))_{i,j} \frac{\partial^2}{\partial x_i \, \partial x_j} + \sum_{i=1}^n b_i(x) \frac{\partial}{\partial x_i}.$$

Consider the following SPDE with nonlinear Neumann boundary condition:

$$(f,g,h) \begin{cases} \mathrm{d}u(t,x) + [Lu(t,x) + f(t,x,u(t,x),\sigma^*(x)D_x u(t,x))] \, \mathrm{d}t \\ \quad + \sum_{i=1}^d g_i(t,x,u(t,x)) \overleftarrow{\mathrm{d}B_t^i} = 0, \qquad (t,x) \in [0,T] \times G, \\ u(T,x) = l(x), \qquad x \in G, \\ \dfrac{\partial u}{\partial n}(t,x) + h(t,x,u(t,x)) = 0, \qquad x \in \partial G. \end{cases} \tag{14}$$

We now define the notion of stochastic viscosity solution for the SPDE $(f,g,h)$. We are inspired by the work of Buckdahn and Ma [2] and we refer to their paper for a lucid discussion on this topic. We use some of their notation and follow the lines of their proofs to obtain our main result. Indeed, we will use the stochastic flow $\widehat{\eta}(t,x,y) \in \mathcal{C}(\mathbf{F}^B, [0,T] \times \mathbb{R}^n \times \mathbb{R})$, defined as the unique solution of the SDE which, in Stratonowich form, reads as follows:

$$\widehat{\eta}(t,x,y) = y + \sum_{i=1}^d \int_0^t g_i(s,x,\widehat{\eta}(s,x,y)) \circ \mathrm{d}B_s^i,$$

$$= y + \int_0^t \langle g(s,x,\widehat{\eta}(s,x,y)), \circ \mathrm{d}B_s \rangle, \qquad t \geq 0. \tag{15}$$

Under the assumption (H$_5$) the mapping $y \mapsto \widehat{\eta}(t,x,y)$ defines a diffeomorphism for all $(t,x)$, $\mathbb{P}$-almost surely see Protter [9]. Denote the $y$-inverse of $\widehat{\eta}(t,x,y)$ by $\widehat{\varepsilon}(t,x,y)$. Then, since $\widehat{\varepsilon}(t,x,\widehat{\eta}(t,x,y)) = y$, one can show that (cf. Buckdahn and Ma [2])

$$\widehat{\varepsilon}(t,x,y) = y - \int_0^t \langle D_y \widehat{\varepsilon}(s,x,y), g(s,x,y) \circ \mathrm{d}B_s \rangle, \tag{16}$$

where the stochastic integrals have to be interpreted in Stratonowich sense. Now let us introduce the process $\eta \in \mathcal{C}(\mathbf{F}^B, [0,T] \times \mathbb{R}^n \times \mathbb{R})$ as the solution to the equation

$$\eta(t,x,y) = y + \int_t^T \langle g(s,x,\eta(s,x,y)), \circ \overleftarrow{\mathrm{d}B}_s \rangle, \qquad 0 \leq t \leq T. \tag{17}$$

We note that due to the direction of the Itô integral, equation (17) should be viewed as going from $T$ to $t$ (i.e. $y$ should be understood as the initial value). Then $y \longmapsto \eta(s,x,y)$ will have the same regularity properties as those of $y \longmapsto \widehat{\eta}(s,x,y)$ for all $(s,x) \in [t,T] \times \mathbb{R}^n$, $\mathbb{P}$-almost surely. Hence if we denote by $\varepsilon$ its $y$-inverse, we obtain

$$\varepsilon(t,x,y) = y - \int_t^T \langle D_y\varepsilon(s,x,y), g(s,x,y) \circ \overleftarrow{\mathrm{d}B}_s \rangle. \tag{18}$$

To simplify the notation we write:

$$\mathrm{A}_{f,g}(\varphi(t,x)) = -L\varphi(t,x) - f(t,x,\varphi(t,x), \sigma^*(x)D_x\varphi(t,x)) + \tfrac{1}{2}\langle g, D_y g \rangle(t,x,\varphi(t,x)).$$

We now introduce the notion of a stochastic viscosity solution of the SPDE $(f,g,h)$ as follows.

***Definition 3.1.*** *A random field* $u \in \mathcal{C}(\mathbf{F}^B, [0,T] \times \overline{\mathrm{G}})$ *is called a stochastic viscosity subsolution of the SPDE$(f,g,h)$ if* $u(T,x) \leq l(x)$, *for all* $x \in \overline{\mathrm{G}}$, *and if for any stopping time* $\tau \in \mathcal{M}_{0,T}^B$, *any state variable* $\xi \in L^0(\mathcal{F}_\tau^B, [0,T] \times \mathrm{G})$, *and any random field* $\varphi \in \mathcal{C}^{1,2}(\mathcal{F}_\tau^B, [0,T] \times \mathbb{R}^n)$, *with the property that for $\mathbb{P}$-almost all $\omega \in \{0 < \tau < T\}$ the inequality*

$$u(t,\omega,x) - \eta(t,\omega,x,\varphi(t,x)) \leq 0 = u(\tau(\omega),\xi(\omega)) - \eta(\tau(\omega),\xi(\omega),\varphi(\tau(\omega),\xi(\omega)))$$

*is fulfilled for all $(t,x)$ in some neighbourhood $\mathcal{V}(\omega,\tau(\omega),\xi(\omega))$ of $(\tau(\omega),\xi(\omega))$, the following conditions are satisfied:*

(a) *On the event $\{0 < \tau < T\}$ the inequality*

$$\mathrm{A}_{f,g}(\psi(\tau,\xi)) - D_y\eta(\tau,\xi,\varphi(\tau,\xi))D_t\varphi(\tau,\xi) \leq 0 \tag{19}$$

*holds $\mathbb{P}$-almost surely.*

(b) *On the event $\{0 < \tau < T\} \cap \{\xi \in \partial \mathrm{G}\}$ the inequality*

$$\min\left[\mathrm{A}_{f,g}(\psi(\tau,\xi)) - D_y\eta(\tau,\xi,\varphi(\tau,\xi))D_t\varphi(\tau,\xi),\right.$$
$$\left. -\frac{\partial \psi}{\partial n}(\tau,\xi) - h(\tau,\xi,\psi(\tau,\xi))\right] \leq 0 \tag{20}$$

*holds $\mathbb{P}$-almost surely with $\psi(t,x) \triangleq \eta(t,x,\varphi(t,x))$.*

A random field $u \in \mathcal{C}(\mathbf{F}^B, [0,T] \times \overline{\mathrm{G}})$ is called a *stochastic viscosity supersolution* of the SPDE $(f,g,h)$ if $u(T,x) \geq l(x)$, for all $x \in \overline{\mathrm{G}}$, and if for any stopping time $\tau \in \mathcal{M}^B_{0,T}$, any state variable $\xi \in L^0(\mathcal{F}^B_\tau, [0,T] \times \mathrm{G})$, and any random field $\varphi \in \mathcal{C}^{1,2}(\mathcal{F}^B_\tau, [0,T] \times \mathbb{R}^n)$, with the property that for $\mathbb{P}$-almost all $\omega \in \{0 < \tau < T\}$ the inequality

$$u(t,\omega,x) - \eta(t,\omega,x,\varphi(t,x)) \geq 0 = u(\tau(\omega), \xi(\omega)) - \eta(\tau(\omega), \xi(\omega), \varphi(\tau(\omega), \xi(\omega)))$$

is fulfilled for all $(t,x)$ in some neighbourhood $\mathcal{V}(\omega, \tau(\omega), \xi(\omega))$ of $(\tau(\omega), \xi(\omega))$, the following conditions are satisfied:

(a) on the event $\{0 < \tau < T\}$ the inequality

$$\mathrm{A}_{f,g}(\psi(\tau,\xi)) - D_y \eta(\tau,\xi,\varphi(\tau,\xi)) D_t \varphi(\tau,\xi) \geq 0 \tag{21}$$

holds $\mathbb{P}$-almost surely;

(b) on the event $\{0 < \tau < T\} \cap \{\xi \in \partial \mathrm{G}\}$ the inequality

$$\max\bigg[\mathrm{A}_{f,g}(\psi(\tau,\xi)) - D_y \eta(\tau,\xi,\varphi(\tau,\xi)) D_t \varphi(\tau,\xi), \\ -\frac{\partial \psi}{\partial n}(\tau,\xi) - h(\tau,\xi,\psi(\tau,\xi))\bigg] \geq 0 \tag{22}$$

holds $\mathbb{P}$-almost surely with $\psi(t,x) \triangleq \eta(t,x,\varphi(t,x))$.

Finally, a random field $u \in \mathcal{C}(\mathbf{F}^B, [0,T] \times \overline{\mathrm{G}})$ is called a *stochastic viscosity solution* of the SPDE $(f,g,h)$ if it is both a stochastic viscosity subsolution and a supersolution.

**Remark 3.1.** Observe that if $f, h$ are deterministic and $g \equiv 0$ then Definition 3.1 coincides with the deterministic case (cf. [8]).

Now let us recall a notion of random viscosity solution which will be a bridge linking the stochastic viscosity solution and its deterministic counterpart.

**Definition 3.2.** A random field $u \in \mathcal{C}(\mathbf{F}^B, [0,T] \times \mathbb{R}^n)$ is called an *$\omega$-wise viscosity solution* if for $\mathbb{P}$-almost all $\omega \in \Omega$, $u(\omega, \cdot, \cdot)$ is a (deterministic) viscosity solution of the SPDE $(f, 0, h)$.

Next we introduce the Doss–Sussman transformation. It enables us to convert an SPDE of the form $(f,g,h)$ to an ordinary differential equation of the form $(\widetilde{f}, 0, \widetilde{h})$, where $\widetilde{f}$ and $\widetilde{h}$ are certain well-defined random fields, which are defined in terms of $f, g, h$.

**Proposition 3.1.** *Assume* $(\mathrm{H}_1)$–$(\mathrm{H}_5)$. *A random field $u$ is a stochastic viscosity solution to the SPDE $(f,g,h)$ if and only if $v(\cdot,\cdot) = \varepsilon(\cdot,\cdot,u(\cdot,\cdot))$ is a stochastic viscosity solution to the SPDE $(\widetilde{f}, 0, \widetilde{h})$, where $(\widetilde{f}, \widetilde{h})$ are two coefficients that will be made precise later (see* (24) *and* (26) *below).*

**Remark 3.2.** Let us recall that under assumption ($H_5$) the random field $\eta$ belongs to $C^{0,2,2}(\mathbf{F}^B, [0,T] \times \mathbb{R}^n \times \mathbb{R})$, and hence that the same is true for $\varepsilon$. Then, considering the transformation $\psi(t,x) = \eta(t,x,\varphi(t,x))$, we obtain

$$D_x\psi = D_x\eta + D_y\eta D_x\varphi,$$

$$D_{xx}\psi = D_{xx}\eta + 2(D_{xy}\eta)(D_x\varphi)^* + (D_{yy}\eta)(D_x\varphi)(D_x\varphi)^* + (D_y\eta)(D_{xx}\varphi).$$

Moreover, since for all $(t,x,y) \in [0,T] \times \mathbb{R}^n \times \mathbb{R}$ the equality $\varepsilon(t,x,\eta(t,x,y)) = y$ holds $\mathbb{P}$-almost surely, we also have

$$D_x\varepsilon + D_y\varepsilon D_x\eta = 0,$$
$$D_y\varepsilon D_y\eta = 1,$$
$$D_{xx}\varepsilon + 2(D_{xy}\varepsilon)(D_x\eta)^* + (D_{yy}\varepsilon)(D_x\eta)(D_x\eta)^* + (D_y\varepsilon)(D_{xx}\eta) = 0,$$
$$(D_{xy}\varepsilon)(D_y\eta) + (D_{yy}\varepsilon)(D_x\eta)(D_y\eta) + (D_y\varepsilon)(D_{xy}\eta) = 0,$$
$$(D_{yy}\varepsilon)(D_y\eta)^2 + (D_y\varepsilon)(D_{yy}\eta) = 0,$$

where all the derivatives of the random field $\varepsilon(\cdot,\cdot,\cdot)$ are evaluated at $(t,x,\eta(t,x,y))$, and all those of $\eta(\cdot,\cdot,\cdot)$ are evaluated at $(t,x,y)$.

**Proof of Proposition 3.1.** We shall only argue for the stochastic subsolution case, as the supersolution part is similar. Therefore, in the present proof we assume that $u \in \mathcal{C}(\mathbf{F}^B, [0,T] \times \overline{G})$ is a stochastic viscosity subsolution of the SPDE $(f,g,h)$. It then follows that $v(\cdot,\cdot) = \varepsilon(\cdot,\cdot,u(\cdot,\cdot))$ belongs to $\mathcal{C}(\mathbf{F}^B, [0,T] \times \overline{G})$. In order to show that $v$ is a stochastic viscosity subsolution of the SPDE $(\widetilde{f},0,\widetilde{h})$, we let $\tau \in \mathcal{M}_{0,T}^B$, $\xi \in L^0(\mathcal{F}_\tau^B, [0,T] \times G)$ and $\varphi \in \mathcal{C}^{1,2}(\mathcal{F}_\tau^B, [0,T] \times \mathbb{R}^n)$ be such that for $\mathbb{P}$-almost all $\omega \in \{0 < \tau < T\}$ the inequality

$$v(t,x) - \varphi(t,x) \leq 0 = v(\tau(\omega),\xi(\omega)) - \varphi(\tau(\omega),\xi(\omega))$$

holds for all $(t,x)$ in some neighbourhood $\mathcal{V}(\omega,\tau(\omega),\xi(\omega))$ of $(\tau(\omega),\xi(\omega))$. Next we put $\psi(t,x) = \eta(t,x,\varphi(t,x))$. Since the mapping $y \mapsto \eta(t,x,y)$ is strictly increasing, for all $(t,x) \in \mathcal{V}(\tau,\xi)$ we have that

$$u(t,x) - \psi(t,x) = \eta(t,x,v(t,x)) - \eta(t,x,\varphi(t,x))$$
$$\leq 0 = \eta(\tau,\xi,v(\tau,\xi)) - \eta(\tau,\xi,\varphi(\tau,\xi))$$
$$= u(\tau,\xi) - \psi(\tau,\xi)$$

holds $\mathbb{P}$-almost surely on $\{0 < \tau < T\}$. Moreover, since $u$ is a stochastic viscosity subsolution of the SPDE $(f,g,h)$, the inequality

$$\mathrm{A}_{f,g}(\psi(\tau,\xi)) - D_y\eta(\tau,\xi,\varphi(\tau,\xi))D_t\varphi(\tau,\xi) \leq 0 \tag{23}$$

holds $\mathbb{P}$-almost everything on the event $\{0 < \tau < T\}$. On the other hand, from Remark 3.2 it follows that

$$L\psi(t,x) = \frac{1}{2}\operatorname{tr}(\sigma(x)\sigma(x)^*D_{xx}\psi(t,x)) + \langle b(x), D_x\psi(t,x)\rangle$$
$$= L_x\eta(t,x,\varphi(t,x)) + D_y\eta(t,x,\varphi(t,x))L\varphi(t,x)$$
$$+ \langle \sigma(x)\sigma^\star(x)D_{xy}\eta(t,x,\varphi(t,x)), D_x\varphi(t,x)\rangle$$
$$+ \frac{1}{2}D_{yy}\eta(t,x,\varphi(t,x))(D_x\varphi)^\star(D_x\varphi),$$

where $L_x$ is the same as the operator $L$ with all the derivatives taken with respect to the second variable $x$. Then if we define the random field $\widetilde{f}$ by

$$\widetilde{f}(t,x,y,z) = \frac{1}{D_y\eta(t,x,y)}\left[f(t,x,\eta(t,x,y),\sigma(x)^\star D_x\eta(t,x,y) + D_y\eta(t,x,y)z)\right.$$
$$- \frac{1}{2}gD_yg(t,x,\eta(t,x,y)) + L_x\eta(t,x,y)$$
$$\left. + \langle \sigma(x)^\star D_{x,y}\eta(t,x,y), z\rangle + \frac{1}{2}D_{yy}\eta(t,x,y)|z|^2\right], \quad (24)$$

we obtain

$$D_y\varepsilon(t,x,\psi(t,x))\mathrm{A}_{f,g}(\psi(t,x)) = \mathrm{A}_{\widetilde{f},0}(\varphi(t,x)).$$

Here we have used the fact that the equality

$$D_y\eta(t,x,\varphi(t,x)) = \frac{1}{D_y\varepsilon(t,x,\psi(t,x))} \qquad \forall (t,x)$$

holds $\mathbb{P}$-almost everywhere for all $(t,x) \in [0,T] \times \mathrm{G}$. Consequently, (23) becomes

$$\mathrm{A}_{\widetilde{f},0}(\varphi(\tau,\xi)) - D_t\varphi(\tau,\xi) \leq 0, \quad (25)$$

and hence the Doss–Sussman transformation converts an SPDE of the form $(f,g,h)$ to one of the form $(f,0,h)$, provided a similar transformation for the random field $h$ also works for the inequality in (22). This establishes part (a) of Definition 3.1.

In order to establish part (b) in Definition 3.1 we notice that for all $(t,x) \in [0,T] \times \partial\mathrm{G}$ the following string of equalities holds:

$$\frac{\partial\psi}{\partial n}(t,x) = \langle D_x\psi(t,x), \nabla\phi(x)\rangle$$
$$= \langle D_x\eta(t,x,\varphi(t,x)), \nabla\phi(x)\rangle + D_y\eta(t,x,\varphi(t,x))\langle D_x\varphi(t,x), \nabla\phi(x)\rangle$$
$$= \langle D_x\eta(t,x,\varphi(t,x)), \nabla\phi(x)\rangle + D_y\eta(t,x,\varphi(t,x))\frac{\partial\varphi}{\partial n}(t,x).$$

Hence,

$$\frac{\partial \psi}{\partial n}(\tau,\xi) + h(\tau,\xi,\psi(\tau,\xi)) = D_y\eta(\tau,\xi,\varphi(\tau,\xi))\frac{\partial \varphi}{\partial n}(\tau,\xi) + \langle D_x\eta(\tau,\xi,\varphi(\tau,\xi)), \nabla\phi(x)\rangle$$
$$+ h(\tau,\xi,\eta(\tau,\xi,\varphi(\tau,\xi)))$$
$$= D_y\eta(\tau,\xi,\varphi(\tau,\xi))\left(\frac{\partial \varphi}{\partial n}(\tau,\xi) + \widetilde{h}(\tau,\xi,\varphi(\tau,\xi))\right),$$

where

$$\widetilde{h}(t,x,y) = \frac{1}{D_y\eta(t,x,y)}(h(t,x,\eta(t,x,y)) + \langle D_x\eta(t,x,y), \nabla\phi(x)\rangle). \qquad (26)$$

Since $D_y\eta(t,x,y) > 0$, we obtain, $\mathbb{P}$-almost surely on the event $\{0 < \tau < T\} \cap \{\xi \in \partial G\}$, the inequality

$$\min\left[A_{\widetilde{f},0}(\varphi(\tau,\xi)) - D_t\varphi(\tau,\xi), -\frac{\partial \varphi}{\partial n}(\tau,\xi) - \widetilde{h}(\tau,\xi,\psi(\tau,\xi))\right] \leq 0. \qquad (27)$$

Combining (25) and the inequality in (27), we obtain that the random field $v$ is a stochastic viscosity subsolution of the SPDE $(\widetilde{f}, 0, \widetilde{h})$, which concludes the proof of Proposition 3.1. $\square$

## 4. Generalized backward doubly SDEs and viscosity solution of SPDEs with nonlinear Neumann boundary conditions

The main objective of this section is to show how a semi-linear SPDE with coefficients $(f, g, h)$ is related to equation (1) introduced in Section 1.

### 4.1. Reflected diffusions

In this section we recall some known results on reflected diffusions. First, since $\sigma$ and $b$ satisfy condition (H$_3$), it follows from Lions and Sznitman [5] that for each $x \in \overline{G}$ there exists a unique pair of progressively measurable continuous processes $(X^x, k^x)$ with values in $\overline{G} \times \mathbb{R}_+$ such that

$$t \mapsto k_t^x \text{ is increasing},$$
$$X_t^x = x + \int_0^t b(X_r^x)\,dr + \int_0^t \sigma(X_r^x)\,dW_r + \int_0^t \nabla\phi(X_r^x)\,dk_r^x, \qquad \text{for } t \in [0,T],$$
$$k_t^x = \int_0^t I_{\{X_s^x \in \partial G\}}\,dk_s^x,$$

where the stochastic integral is the standard Itô integral, and the probability space (and its filtration) is the one on which the Brownian motion $W$ is defined. We refer to Pardoux and Zhang [8]: Propositions 3.1 and 3.2, for the following regularity results.

**Proposition 4.1.** *There exist a constant $C > 0$ such that for all $x$, $x' \in \overline{G}$ the following inequality holds:*

$$\mathbb{E}\left[\sup_{0 \leq t \leq T} |X_t^x - X_t^{x'}|^4\right] \leq C|x - x'|^4.$$

**Proposition 4.2.** *For each $T > 0$, there exists a constant $C_T$ such that for all $x$, $x' \in G$,*

$$\mathbb{E}\left(\sup_{0 \leq t \leq T} |k_t^x - k_t^{x'}|^4\right) \leq C_T |x - x'|^4.$$

*Moreover, for all $p \geq 1$, there exists a constant $C_p$ such that for all $(t, x) \in \mathbb{R}_+ \times \overline{G}$,*

$$\mathbb{E}(|k_t^x|^p) \leq C_p(1 + t^p)$$

*and for each $\mu$, $t > 0$, there exists a constant $C(\mu, t)$ such that for all $x \in \overline{G}$,*

$$\mathbb{E}(e^{\mu k_t^x}) \leq C(\mu, t).$$

### 4.2. Viscosity solutions

Consider, for every $(t, x) \in [0, T] \times \overline{G}$, the process $s \mapsto (X_s^{t,x}, k_s^{t,x})$, $s \in [0, T]$, as the unique solution of the equation

$$X_s^{t,x} = x + \int_t^{s \vee t} b(X_r^{t,x})\,dr + \int_t^{s \vee t} \sigma(X_r^{t,x})\,dW_r + \int_t^{s \vee t} \nabla\phi(X_r^{t,x})\,dk_r^{t,x}.$$

The main topic of this section will be a study of properties of the solution $(Y_s^{t,x}, Z_s^{t,x})$, $(t, x) \in [0, T] \times \overline{G}$, to the BDSDE

$$Y_s^{t,x} = l(X_T^{t,x}) + \int_{s \vee t}^T f(r, X_r^{t,x}, Y_r^{t,x}, Z_r^{t,x})\,dr + \sum_{i=1}^d \int_{s \vee t}^T g_i(r, X_r^{t,x}, Y_r^{t,x}) \overleftarrow{dB_r^i}$$

$$+ \int_{s \vee t}^T h(r, X_r^{t,x}, Y_r^{t,x})\,dk_r^{t,x} - \int_s^T \langle Z_r^{t,x}, dW_r\rangle, \qquad 0 \leq s \leq T, \qquad (28)$$

where the coefficients $l$, $f$, $g$ and $h$ satisfy the hypotheses $(H_1')$, $(H_2)$, $(H_4)$ and $(H_5)$.

**Proposition 4.3.** *Let the ordered pair $(Y_s^{t,x}, Z_s^{t,x})$ be a solution to the BSDE in (28). Then the random field $(s, t, x) \to Y_s^{t,x}$, $(s, t, x) \in [0, T] \times [0, T] \times \overline{G}$ is almost surely continuous.*

**Proof.** Let $(t, x)$ and $(t', x')$ be elements of $[0, T] \times \overline{G}$. It follows from Remark 2.2 that, for $0 \leq s \leq T$,

$$|Y_s^{t,x} - Y_s^{t',x'}|^2$$
$$\leq C \mathbb{E}^{\mathcal{F}_s} \bigg( |e^{\mu A_T}[l(X_T^{t,x}) - l(X_T^{t',x'})]|^2 + \int_0^T e^{\mu A_s} |X_s^{t,x} - X_s^{t',x'}|^2 \, ds$$
$$+ \int_0^T e^{\mu A_s} |X_s^{t,x} - X_s^{t',x'}|^2 \, dk_s^{t',x'}$$
$$+ \sup_{0 \leq s \leq T} (1 + |X_s^{t,x}|^2 + |Y_s^{t,x}|^2) e^{\mu A_T} |k^{t,x} - k^{t',x'}|_T \bigg).$$

The result follows from standard arguments using Propositions 4.1, 4.2, 2.4 and the continuity of the function $l$. □

Next we recall a generalized version of the Itô–Ventzell formula; the proof is analogous to the corresponding one in Buckdahn and Ma [2].

**Theorem 4.4.** *Suppose that $M \in \mathcal{C}^{0,2}(F, [0, T] \times \mathbb{R}^n)$ is a semimartingale in the sense that, for every spatial parameter $x \in \mathbb{R}^n$, the process $t \mapsto M(t, x)$, $t \in [0, T]$, is of the form*

$$M(t, x) = M(0, x) + \int_0^t G(s, x) \, ds + \int_0^t \langle H(s, x), \overleftarrow{dB_s} \rangle + \int_0^t \langle K(s, x), dW_s \rangle,$$

*where $G \in \mathcal{C}^{0,2}(\mathbf{F}^B, [0, T] \times \mathbb{R}^n)$ and $H \in \mathcal{C}^{0,2}(\mathbf{F}^B, [0, T] \times \mathbb{R}^n; \mathbb{R}^d)$, and the process $K$ belongs to $\mathcal{C}^{0,2}(\mathbf{F}^W, [0, T] \times \mathbb{R}^n; \mathbb{R}^d)$. Let $\alpha \in \mathcal{C}(F, [0, T]; \mathbb{R}^n)$ be a process of the form*

$$\alpha_t = \alpha_0 + \int_0^t \beta_s \, dk_s + \int_0^t \gamma_s \overleftarrow{dB_s} + \int_0^t \delta_s \, dW_s,$$

*where $\beta \in \mathcal{K}^2(F, [0, T]; \mathbb{R}^n)$, $\gamma \in \mathcal{M}^2(F, [0, T]; \mathbb{R}^{n \times d})$ and $\delta \in \mathcal{M}^2(F, [0, T]; \mathbb{R}^{n \times d})$. Then the following equality holds $\mathbb{P}$-almost surely for all $0 \leq t \leq T$:*

$$M(t, \alpha_t) = M(0, \alpha_0) + \int_0^t G(s, \alpha_s) \, ds + \int_0^t \langle H(s, \alpha_s), \overleftarrow{dB_s} \rangle + \int_0^t \langle K(s, \alpha_s), dW_s \rangle$$
$$+ \int_0^t \langle D_x M(s, \alpha_s), \beta_s \, dk_s \rangle + \int_0^t \langle D_x M(s, \alpha_s), \gamma_s \overleftarrow{dB_s} \rangle + \int_0^t \langle D_x M(s, \alpha_s), \delta_s \, dW_s \rangle$$
$$+ \frac{1}{2} \int_0^t \mathrm{tr}(D_{xx} M(s, \alpha_s) \delta_s \delta_s^*) \, ds - \frac{1}{2} \int_0^t \mathrm{tr}(D_{xx} M(s, \alpha_s) \gamma_s \gamma_s^*) \, ds$$
$$+ \int_0^t \mathrm{tr}(D_x K(s, \alpha_s) \delta_s^*) \, ds - \int_0^t \mathrm{tr}(D_x H(s, \alpha_s) \gamma_s^*) \, ds.$$

### 4.3. Existence of stochastic viscosity solutions

In this section we apply the results of the previous sections to prove the existence of stochastic viscosity solutions to a quasi-linear SPDE with Neumann boundary conditions. To this end, we need the following result which is proved in Buckdahn and Ma [2]:

**Proposition 4.5.** *Assume* $(H_5)$. *Let $\eta$ be the unique solution to SDE* (15) *and $\varepsilon$ be the $y$-inverse of $\eta$. Then there exists a constant $C > 0$, depending only on the bound of $g$ and its partial derivatives, such that for $\zeta = \eta$, and $\zeta = \varepsilon$, the following inequalities hold $\mathbb{P}$-almost surely for all $(t, x, y) \in [0, T] \times \mathbb{R}^n \times \mathbb{R}$:*

$$|\zeta(t,x,y)| \leq |y| + C|B_t|,$$

$$|D_x\zeta|, \ |D_y\zeta|, \ |D_{xx}\zeta|, \ |D_{xy}\zeta|, \ |D_{yy}\zeta| \leq C\exp\{C|B_t|\}.$$

*Here all the derivatives are evaluated at $(t, x, y)$.*

Next, for $t \in [0, T]$ and $x \in \overline{G}$, let us define the following processes:

$$U_s^{t,x} = \varepsilon(s, X_s^{t,x}, Y_s^{t,x}), \qquad 0 \leq t \leq s \leq T,$$

$$V_s^{t,x} = D_y\varepsilon(s, X_s^{t,x}, Y_s^{t,x})Z_s^{t,x} + \sigma^*(X_s^{t,x})D_x\varepsilon(s, X_s^{t,x}, Y_s^{t,x}), \qquad 0 \leq t \leq s \leq T.$$

Then from Proposition 4.5 we obtain

$$((U_s^{t,x}, V_s^{t,x}), (s, x) \in [0, T] \times \overline{G}) \in \mathcal{S}^2(F; [0, T]; \mathbb{R}) \times \mathcal{M}^2(F; [0, T]; \mathbb{R}^d).$$

**Theorem 4.6.** *For each $(t, x) \in [0, T] \times \overline{G}$, the process $(U_s^{t,x}, V_s^{t,x}, \ t \leq s \leq T)$ is the unique solution to the following generalized BSDE:*

$$U_s^{t,x} = l(X_T^{t,x}) + \int_s^T \widetilde{f}(t, X_r^{t,x}, U_r^{t,x}, V_r^{t,x})\,\mathrm{d}r$$

$$+ \int_s^T \widetilde{h}(r, X_r^{t,x}, U_r^{t,x},)\,\mathrm{d}k_r^{t,x} - \int_s^T V_r^{t,x}\,\mathrm{d}W_r, \tag{29}$$

*where $\widetilde{f}$ and $\widetilde{h}$ are given by* (24) *and* (26).

**Proof.** For brevity we write $X$, $Y$, $U$, $V$, $Z$ and $k$ instead of $X^{t,x}$, $Y^{t,x}$, $U^{t,x}$, $V^{t,x}$, $Z^{t,x}$ and $k^{t,x}$, respectively. Then the mapping $(X, Y, Z) \mapsto (X, U, V)$ is one-to-one with inverse transformation

$$Y_s = \eta(s, X_s, U_s), \qquad Z_s = D_y\eta(s, X_s, U_s)V_s + \sigma^*(X_s)D_x\eta(s, X_s, U_s).$$

Then the uniqueness of solutions to the equation in Theorem 4.6 follows from that of the generalized BDSDE in (28). As a consequence we only need to show that $(U, V)$ is a solution of the generalized BSDE in Theorem 4.6. Indeed, using the Itô–Ventzell formula, we obtain

$$U_s = l(X_T) - \int_s^T \langle D_x\varepsilon(r, X_r, Y_r), b(X_r)\rangle \, dr - \int_s^T \langle D_x\varepsilon(r, X_r, Y_r), \sigma(X_r) \, dW_r\rangle$$
$$- \int_s^T \langle D_x\varepsilon(r, X_r, Y_r), \nabla\phi(X_r)\rangle \, dk_r - \tfrac{1}{2}\int_s^T \mathrm{tr}\{\sigma(X_r)\sigma^*(X_r)D_{xx}\varepsilon(r, X_r, Y_r)\} \, dr$$
$$+ \int_s^T D_y\varepsilon(r, X_r, Y_r)f(r, X_r, Y_r, Z_r) \, dr + \int_s^T D_y\varepsilon(r, X_r, Y_r)h(r, X_r, Y_r) \, dk_r$$
$$- \int_s^T \langle D_y\varepsilon(r, X_r, Y_r)Z_r, dW_r\rangle - \tfrac{1}{2}\int_s^T D_{yy}\varepsilon(r, X_r, Y_r)\|Z_r\|^2 \, dr$$
$$- \int_s^T \langle \sigma^*(X_r)D_{xy}\varepsilon(r, X_r, Y_r), Z_r\rangle \, dr - \tfrac{1}{2}\int_s^T D_y\varepsilon(r, X_r, Y_r)\langle g, D_y g\rangle(r, X_r, Y_r) \, dr$$
$$= l(X_T) + \int_s^T \mathcal{F}(r, X_r, Y_r, Z_r) \, dr + \int_s^T \mathcal{H}(r, X_r, Y_r) \, dk_r - \int_s^T \langle V_r, dW_r\rangle, \quad (30)$$

where

$$\mathcal{F}(s, x, y, z) \triangleq -\langle D_x\varepsilon, b(x)\rangle + (D_y\varepsilon)f(s, x, y, z) - \tfrac{1}{2}(D_{yy}\varepsilon)|z|^2$$
$$- \tfrac{1}{2}\mathrm{tr}\{\sigma(x)\sigma(x)^* D_{xx}\varepsilon\} - \langle \sigma^*(x)D_{xy}\varepsilon, z\rangle - \tfrac{1}{2}D_y\varepsilon\langle g, D_y g\rangle(s, x, y) \quad (31)$$

and

$$\mathcal{H}(s, x, y) \triangleq -\langle D_x\varepsilon, \nabla\phi(x)\rangle + (D_y\varepsilon)h(s, x, y). \quad (32)$$

From (30), (31) and (32) it follows that it suffices to show that

$$\mathcal{F}(s, X_s, Y_s, Z_s) = \widetilde{f}(s, X_s, U_s, V_s) \qquad \forall s \in [0, T], \mathbb{P}\text{-a.s.} \quad (33)$$

and

$$\mathcal{H}(s, X_s, Y_s) = \widetilde{h}(s, X_s, U_s) \qquad \forall s \in [0, T], \mathbb{P}\text{-a.s.} \quad (34)$$

To this end, if we write $\sigma(X_s) = \sigma_s$ and $b(X_s) = b_s$, Remark 3.2 entails the following equalities:

$$\langle D_x\varepsilon(s, X_s, Y_s), b_s\rangle = -D_y\varepsilon(s, X_s, Y_s)\langle D_x\eta(s, X_s, U_s), b_s\rangle,$$
$$(D_y\varepsilon)f(s, X_s, Y_s, Z_s) = (D_y\varepsilon)f(s, X_s, \eta, D_y\eta V_s + \sigma_s^*(D_x\eta)),$$
$$\langle \sigma_s^*(D_{x,y}\varepsilon), Z_s\rangle = (D_y\eta)\langle \sigma_s^*(D_{x,y}\varepsilon), V_s\rangle + \langle \sigma_s^*(D_{x,y}\varepsilon), \sigma_s^*(D_x\eta)\rangle,$$
$$-\tfrac{1}{2}(D_{yy}\varepsilon)|Z_s|^2 = \tfrac{1}{2}(D_y\varepsilon)(D_{yy}\eta)|V_s|^2 + (D_y\varepsilon)^2(D_{yy}\eta)\langle V_s, \sigma_s^*(D_x\eta)\rangle$$
$$+ \tfrac{1}{2}(D_{yy}\eta)(D_y\varepsilon)|\sigma_s^*(D_x\eta)(D_y\varepsilon)|^2. \quad (35)$$

*Generalized backward doubly stochastic differential equations* 445Hence, from the equalities in (35) we obtain

$$\begin{aligned}\mathcal{F}(s,X_s,Y_s,Z_s) &= D_y\varepsilon[\langle D_x\eta,b_s\rangle + \tfrac{1}{2}(D_{yy}\eta)|V_s|^2 \\
&\quad + f(s,X_s,\eta,D_y\eta V_s + \sigma_s^*(D_x\eta)) - \tfrac{1}{2}\langle g,D_y g\rangle(s,X_s,\eta)] \\
&\quad + \langle V_s, \sigma_s^*[D_x\eta(D_y\varepsilon)^2(D_{yy}\eta) - D_y\eta D_{xy}\varepsilon]\rangle, \\
&\quad + [\tfrac{1}{2}(D_{yy}\eta)(D_y\varepsilon)|\sigma_s^*(D_x\eta)(D_y\varepsilon)|^2 \\
&\quad - \tfrac{1}{2}\operatorname{tr}\{\sigma_s\sigma_s^*D_{xx}\varepsilon\} - \langle\sigma_s^*(D_{x,y}\varepsilon),\sigma_s^*(D_x\eta)\rangle],\end{aligned} \qquad (36)$$

where all the derivatives of the random field $\varepsilon(\cdot,\cdot,\cdot)$ are to be evaluated at the point $(s,x,\eta(s,x,y))$, and all those of $\eta(\cdot,\cdot,\cdot)$ at $(s,x,y)$.

Now from Remark 3.2, we have

$$\operatorname{tr}\{\sigma_s\sigma_s^*D_{xx}\varepsilon\} = -2\langle\sigma_s^*D_{xy}\varepsilon,\sigma_s^*D_x\eta\rangle + (D_y\varepsilon)D_{yy}\eta|\sigma_s^*D_x\eta D_y\varepsilon|^2 \\
- (D_y\varepsilon)\operatorname{tr}\{(\sigma_s\sigma_s^*D_{xx}\eta)\} \qquad (37)$$

and

$$D_{xy}\varepsilon D_y\eta - D_x\eta(D_y\varepsilon)^2(D_{yy}\eta) = -D_y\varepsilon D_{xy}\eta. \qquad (38)$$

The equalities in (37) and (38), together with $D_y\varepsilon(s,X_s,Y_s) = (D_y\eta)^{-1}(s,X_s,U_s)$, imply

$$\begin{aligned}\mathcal{F}(s,X_s,Y_s,Z_s) &= D_y\varepsilon[\langle D_x\eta,b_s\rangle + \tfrac{1}{2}(D_{yy}\eta)|V_s|^2 \\
&\quad + f(s,X_s,\eta,D_y\eta V_s + \sigma_s^*(D_x\eta)) - \tfrac{1}{2}\langle g,D_y g\rangle(s,X_s,\eta)] \\
&\quad + \tfrac{1}{2}(D_y\varepsilon)\operatorname{tr}\{\sigma_s\sigma_s^*D_{xx}\eta\} + (D_y\varepsilon)\langle V_s,\sigma_s^*D_{xy}\eta\rangle.\end{aligned} \qquad (39)$$

Since the expressions in (36) and (39) are equal, this shows the equality in (33).

The next argument shows the equality in (34):

$$\begin{aligned}\mathcal{H}(s,X_s,Y_s) &= -\langle D_x\varepsilon(r,X_r,Y_r),\nabla\phi(X_r)\rangle + D_y\varepsilon(r,X_r,Y_r)h(r,X_r,Y_r) \\
&= D_y\varepsilon(s,X_s,Y_s)(\langle D_x\eta(s,X_s,U_s),\nabla\phi(X_r)\rangle + h(s,X_s,\eta(s,X_s,U_s))) \\
&= \frac{1}{D_y\eta(s,X_s,U_s)}[h(s,X_s,\eta(s,X_s,U_s)) + \langle D_x\eta(s,X_s,U_s),\nabla\phi(X_s)\rangle] \\
&= \widetilde{h}(s,X_s,U_s).\end{aligned}$$

This completes the proof of Theorem 4.6. $\square$

To conclude this section, we give our main result. Define for each $(t,x) \in [0,T] \times \overline{G}$ the (random) fields $u$ and $v$ by $u(t,x) = Y_t$ and $v(t,x) = U_t$, where $(Y,Z)$ and $(U,V)$ are the solutions to the BSDEs (28) and (29), respectively. Then we have

$$u(t,\omega,x) = \eta(\omega,t,x,v(t,\omega,x)), \qquad v(t,\omega,x) = \varepsilon(\omega,t,x,u(t,\omega,x)). \qquad (40)$$

**Theorem 4.7.** *Under the above assumptions, the random field u is a stochastic viscosity solution to the SPDE $(f, g, h)$.*

***Remark 4.1.*** From the results in Proposition 4.5 we see that, in order to prove Theorem 4.7, we only need to show that the random field $v$ is a stochastic viscosity solution to the SPDE $(\widetilde{f}, 0, \widetilde{h})$.

**Proof of Theorem 4.7.** From Proposition 4.3 it follows that the mapping $(s, t, x) \mapsto Y_s^{t,x}$ is continuous, for all $(s, t, x) \in [0, T]^2 \times \overline{G}$. It follows that $u(t, x) = Y_t^{t,x}$ is continuous as well, and, in particular, it is jointly measurable.

Since $Y_s^{t,x}$ is $\mathcal{F}_{t,s}^W \otimes \mathcal{F}_{t,T}^B$-measurable, it follows that $Y_t^{t,x}$ is $\mathcal{F}_{t,T}^B$-measurable. Consequently, $u(t, x)$ is $\mathcal{F}_{t,T}^B$-measurable and so it is independent of $\omega_1 \in \Omega_1$ (see the notation in Section 2.1). Therefore, we obtain $u \in \mathcal{C}(\mathbf{F}^B, [0, T] \times \overline{G})$, which by (40) implies that $v \in \mathcal{C}(\mathbf{F}^B, [0, T] \times \overline{G})$. However, from Definition 3.2, we see that an $\mathbf{F}^B$-progressively measurable $\omega$-wise viscosity solution is automatically a stochastic viscosity solution. Therefore it suffices to show that $v$ is an $\omega$-wise viscosity solution to the SPDE $(\widetilde{f}, 0, \widetilde{h})$. To this end, we denote, for a fixed $\omega_2 \in \Omega_2$,

$$\overline{U}_s^{\omega_2}(\omega_1) = U_s(\omega_2, \omega_1), \qquad \overline{V}_s^{\omega_2}(\omega_2, \omega_1) = V_s(\omega_2, \omega_1).$$

Since the pair $(\overline{U}^{\omega_2}, \overline{V}^{\omega_2})$ is the unique solution of the generalized BSDE with coefficients $(\widetilde{f}(\omega_2, \cdot, \cdot, \cdot), \widetilde{h}(\omega_2, \cdot, \cdot))$, it follows from Pardoux and Zhang [8] that $\overline{v}(\omega_2, t, x) \triangleq \overline{U}_t^{\omega_2}$ is a viscosity solution to $\text{PDE}(\widetilde{f}(\omega_2, \cdot, \cdot, \cdot), \widetilde{h}(\omega_2, \cdot, \cdot))$ with Neumann boundary condition. By Blumenthal's 0–1 law we have $\mathbb{P}(\overline{U}_t^{\omega_2} = U_t(\omega_2, \omega_1)) = 1$ and, hence, the equality $\overline{v}(t, x) = v(t, x)$ holds $\mathbb{P}_1$-almost surely for all $(t, x) \in [0, T] \times \overline{G}$. Consequently, for every fixed $\omega_2$ the function $v \in \mathcal{C}(\mathbf{F}^B, [0, T] \times \overline{G})$ is a viscosity solution to the SPDE $(\widetilde{f}(\omega_2, \cdot, \cdot, \cdot), 0, \widetilde{h}(\omega, \cdot, \cdot))$. Hence, by definition it is an $\omega$-wise viscosity solution, and thus the result follows from Remark 4.1. This completes the proof. □

## Acknowledgements

The first author wishes to thank the Centre de Recerca Matemàtica for support and hospitality during his stay in July 2005. This work is partially supported by the Fund for Scientific Research Flanders. The third author is grateful to the Moroccan national education department for its financial support. An anonymous referee is acknowledged for his comments, remarks and for a significant improvement to the overall presentation of this paper.